\documentclass{svmult}
  \usepackage{mathptmx}
  \usepackage{helvet}
  \usepackage{courier}
  \usepackage{makeidx}
  \usepackage{graphicx}
  \usepackage{multicol}
  \usepackage{footmisc}
  \usepackage{amsmath} % \boldsymbol
  \usepackage{amssymb} % \mathbb
  \usepackage{color}
  \usepackage{subfig} % subfigures
  \usepackage{tikz}
  \usetikzlibrary{calc,decorations.pathreplacing,intersections,shapes,petri,plotmarks,arrows} % ruins conv plots ch. pot!: positioning
  \usepackage{upgreek} % upepsilon (damn style file)
  \colorlet{verylightshade}{black!4}
  \colorlet{lightshade}{black!8}
  \colorlet{shade}{black!16}
  \colorlet{darkshade}{black!32}
  \colorlet{verydarkshade}{black!50}
  \pgfdeclarehorizontalshading{rainbowh}{100bp}
    {color(0bp)=(blue); color(25bp)=(blue); color(38bp)=(cyan); color(50bp)=(green);
    color(62bp)=(yellow); color(75bp)=(red); color(100bp)=(red)}
  \newcommand{\config}{\theta}
  \newcommand{\constit}{{\Xi}}
  \newcommand{\dom}{\Omega}
    \newcommand{\bnd}{{\partial\dom}}
  \newcommand{    \gammagradi}[3]{#1_{#2;#3}{}}
  \newcommand{\Gin}{{\Gamma_\mathrm{in}}}
  \newcommand{\Gwet}{{\Gamma_t}}
  \newcommand{\refconfig}[1]{{\underline{#1}}}
    \newcommand{\Gwetref}{{\refconfig{\Gamma}}}
  \newcommand{\meass}{s}
  \newcommand{  \ddmu}{\,d\meass}
  \newcommand{\IR}{\mathbb{R}}
  \newcommand{\trac}{\tau}
    
    \newcommand{\reftracj}{\refconfig{\trac}_\mathrm{j}{}}
  \newcommand{  \xinit}{\chi_0}
    \newcommand{\vinit}{\chi_1}
  \newcommand{  \snapshot}[4]{
    \begin{tikzpicture}[shading=rainbowh,scale=0.2,font=\footnotesize]
    \shade (0,0) rectangle node[above] {#1} (20,1);
    \draw (0,1) -- (0,-.5) node[below] {#2};
    \draw (20,1) -- (20,-.5) node[below] {#3};
    \node at (10,3) [anchor=south] {#4};
    \end{tikzpicture}}
\graphicspath{{media/}}
\begin{document}
\title{Finite-element/boundary-element coupling for inflatables: effective contact resolution}
\author{T.M. van Opstal\thanks{
  Department of Mathematical Sciences,
  Norwegian University of Science and Technology,
  7491 Trondheim, Norway,
  \texttt{timo.vanopstal@math.ntnu.no}}}
\maketitle
% \keywords{Boundary element method} % required by Kenji/Yuri?
% 1st paragraph is abstract
\section{Introduction}
A fluid-structure interaction technique for the simulation of inflatable structures is introduced. The presented finite-element/boundary-element (FEBE) technique couples an isogeometric finite element discretization of a flexible shell structure with an isogeometric boundary element discretization of a Stokes fluid. This technique was introduced in~\cite{cm12} for planar problems and extended to 3D in~\cite{3d}. One of the marked advantages of this approach is the contact mechanism, lubrication, which is an inherent attribute of the flow model. Its role in preventing contact was theoretically substantiated, but only demonstrated in the planar setting. In the present work, we demonstrate the effectiveness of the contact mechanism in the isogeometric and three-dimensional setting, discussing various aspects pertaining to the accurate resolution of the traction responsible for contact prevention.
\par
An extensive body of literature already covers the interface of the fields of fluid-structure interaction and contact mechanics. If the problem under consideration requires that the gap between two contact surfaces may vanish, an interface tracking technique may be applied, sometimes locally~\cite{Johnson94,Tezduyar10}, e.g.~\cite{Hsu14a,Kamensky15,Wick13b}. If however, a problem allows for a finite gap to be maintained between the contacting surfaces, Arbitrary Lagrangian-Eulerian~\cite{Hughes81} and Space-time~\cite{Tezduyar92} techniques can been used to compute problems as challenging as the disreefing of parachute clusters~\cite{Takizawa12b,Takizawa11,Takizawa15} and 1000 spheres falling through a tube~\cite{Johnson99}. In this work a finite (but arbitrarily small) gap will likewise be maintained in contact regions.
\par
The target application for the current approach is inflatable structures. These typically undergo large deformations with ubiquitous self-contact during the inflation process, which often starts from a complex, folded initial configuration. Correct simulation is contingent to the resolution of every single contact mode throughout the process. To tackle this type of problem, contact treatment based on the lubrication effect inherent to viscous flow is highly advantageous:
\begin{enumerate}
  \item it avoids explicit contact detection, which is becomes prohibitively expensive and non-robust for complex geometries such as folded inflatables.
  \item accuracy is not impaired by artificial terms added to the model for soft contact treatment, such as in e.g.~\cite{cmame12}.
  \end{enumerate}
Correct simulation of the contact mechanics is contiguous to the accurate resolution of the lubrication tractions, which increase sinularly as the contact gap closes. Thus, after recapulating the mathematical model in~\S\ref{sec:model}, accurate approximation and solution through isogeometric analysis and adaptive quadrature schemes are covered in~\S\ref{sec:approx}. The resulting methodology is demonstrated in a numerical test case, a deflating balloon, in~\S\ref{sec:numex}. Finally, conclusions are drawn in~\S\ref{sec:concl}.
\section{Mathematical model of an inflatable structure}\label{sec:model}
\begin{figure}
  \centering
  \begin{tikzpicture}[scale=1.0,>=latex]
  % Coordinate system
  \begin{scope}[xshift=40mm,yshift=19mm]
    \draw[->] (0, 0) -- (-0.5, -0.2) node [left] {$x_0$};
    \draw[->] (0, 0) -- ( 0.5, -0.2) node [right] {$x_1$};
    \draw[->] (0, 0) -- ( 0.0,  0.5) node [above] {$x_2$};
    \end{scope}
  % Current configuration
  \begin{scope}
    \coordinate (lefthinge) at (-0.6, 0.0);
    \coordinate (rghthinge) at ( 0.6, 0.0);
    % Bulk
    \shadedraw[top color=white, left color=white, right color=lightshade, bottom color=shade, name path=curconfig]
          (lefthinge) to [out= 90, in=-90]
          (-2.0, 1.0) to [out= 90, in=210]
          (-1.8, 1.2) to [out= 30, in=-80]
          (-2.3, 1.8) to [out=100, in=180]
          ( 0.3, 3.0) to [out=  0, in=100]
          ( 2.3, 1.7) to [out=-80, in= 20]
          ( 1.3, 0.5) to [out=200, in= 90]
          (rghthinge) to [out=-90, in=-30]
          ( 0.0, 0.0) to [out=150, in= 90]
          (lefthinge) -- cycle;
    % Fold
    \draw[]
          (-1.8, 1.2) to [out= 30, in=-70]
          (-1.6, 1.6);
    % Inlet, shade
    \draw[verydarkshade, fill=verylightshade] % , densely dashed]
          (rghthinge) to [out=-90, in=-30]
          ( 0.0, 0.0) to [out= 30, in= 90]
          (rghthinge) -- cycle;
    % Inlet, white
    \draw[]
          (rghthinge) to [out=-90, in=-30]
          ( 0.0, 0.0) to [out=210, in=-90]
          (lefthinge);
    % Labels
    \node at ( 0.0, 1.8) {$\dom_t$};
    \draw[->] ( 2.4, 1.1) -- node[pos=0.0, right] {$\Gwet$} ( 2.0, 1.2);
    \coordinate (xcur) at (-1.8, 1.8);
    \end{scope}
  % Reference configuration
  \begin{scope}[xshift=0mm, yshift=0mm]
    \coordinate (lefthinge) at (-0.6, 0.0);
    \coordinate (rghthinge) at ( 0.6, 0.0);
    \coordinate (xref) at (-2.7, 0.6);
    % Phantom DIRTY SOLUTION
    \path[name path=curconfig]
          (lefthinge) to [out= 90, in=-90]
          (-2.0, 1.0) to [out= 90, in=210]
          (-1.8, 1.2) to [out= 30, in=-80]
          (-2.3, 1.8) to [out=100, in=180]
          ( 0.3, 3.0) to [out=  0, in=100]
          ( 2.3, 1.7) to [out=-80, in= 20]
          ( 1.3, 0.5) to [out=200, in= 90]
          (rghthinge);
    \path[name path=refconfig0] (-0.8, 0.4) -- (-3.0, 0.4);
    \path[name path=refconfig1] (-3.0, 0.8) -- ( 0.0, 0.8);
    \path[name path=refconfig2] ( 0.0, 0.8) -- ( 3.0, 0.8);
    \path[name path=refconfig3] ( 3.0, 0.4) -- ( 0.8, 0.4);
    % Bulk
    \shade[top color=verydarkshade, bottom color=verylightshade, semitransparent]
          (lefthinge) --
          (-0.6, 0.2) arc (0:90:2mm) --
          ( 0.8, 0.4) arc (90:180:2mm) --
          (rghthinge) to [out=-90, in=-30]
          ( 0.0, 0.0) to [out=150, in= 90]
          (lefthinge) -- cycle;
    \shade[top color=white, bottom color=verydarkshade, semitransparent]
          (-3.0, 0.4) arc (270:90:2mm) --
          ( 3.0, 0.8) arc (90:-90:2mm) --
          cycle;
    \draw[verydarkshade, name intersections={of=refconfig0 and curconfig, by=x0}]
          (lefthinge) -- (-0.6, 0.2) arc (0:90:2mm) -- (x0);
    \draw[name intersections={of=refconfig1 and curconfig, by=x1}]
          (x0) -- (-3.0, 0.4) arc (270:90:2mm) -- (x1);
    \draw[verydarkshade, name intersections={of=refconfig2 and curconfig, by=x2}]
          (x1) -- (x2);
    \draw[name intersections={of=refconfig3 and curconfig, by=x3}]
          (x2) -- ( 3.0, 0.8) arc (90:-90:2mm) -- (x3);
    \draw[verydarkshade]
          (x3) -- ( 0.8, 0.4) arc (90:180:2mm) -- (rghthinge);
    % Label
    \draw[<-] ( 1.8, 0.45) -- node[pos=1.0, right] {$\Gwetref$} ( 2.2, 0.35);
    \draw[->] (-0.8,-0.1) -- node[pos=0.0, left] {$\Gin$} (-0.4, 0.0);
    \end{scope}
  % Maps
  \begin{scope}
    \draw[->] (xref) to [out=120, in=280] node[midway,above] {$\config_t$} (xcur);
    \draw[fill] (xcur) circle (0.2mm);
    \draw[fill] (xref) circle (0.2mm);
    \end{scope}
  % Additional overlays
  \begin{scope}
    % Inlet
    \draw[darkshade]
          ( 0.0, 0.0) to [out= 30, in= 90]
          (rghthinge);
    \draw[]
          (rghthinge) to [out=-90, in=-30]
          ( 0.0, 0.0) to [out=210, in=-90]
          (lefthinge);
    \end{scope}
  \end{tikzpicture}
  \caption{Schematic geometry of an inflatable structure for the case $d=3$. The reference and current domains are overlayed, as they coincide over $\Gin$.}
  \label{fig:schem}
  \end{figure}
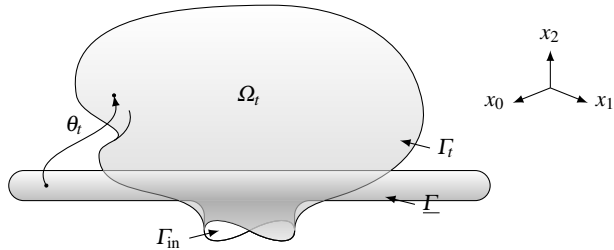

The governing equations of the fluid-structure interaction are elaborated in~\cite{3d}. In this section, this model is recapitulated concisely, and we refer to~\cite{3d} for full details. An inflatable structure occupies a boundary segment $\Gwet$ at time $t$, which, together with a fixed inflow boundary segment $\Gin$ encloses the interior domain $\dom_t$ of the inflatable. This is schematized in fig.~\ref{fig:schem}. The configuration $\config: (0,T)\times\Gwetref\mapsto\Gwet$ maps a material point in the reference manifold $x\in\Gwetref$ to its position on the current manifold $\config(t,x)\in\Gwet$, and is sometimes abbreviated as $\config_t(\cdot):=\config(t,\cdot)$. Throughout, entities related to the reference manifold are underlined.
\par
The flow on both the interior and exterior domains, resp. $\dom_t$ and $\IR^3-\overline{\dom_t}$, are governed by the boundary integral formulation introduced in~\S\ref{sec:fluid}. The motion of the inflatable structure itself is governed by a shell formulation, cf.~\S\ref{sec:structure}. Finally, the coupling of these two subsystems, as well as a condition enforcing compatibility between fluid and structure solutions, is treated in~\S\ref{sec:coupling}.
\subsection{Boundary integral formulation of the fluid}\label{sec:fluid}
The flow interior and exterior to the flexible structure is described by Stokes flow, supplemented by Dirichlet conditions $g$ on the boundary $\bnd_t$ and suitable radiation conditions in the far-field. It will be more convenient to pull back the formulation to the reference configuration, denoting the boundary condition by $\refconfig{g}=g\circ\bnd_t$ and the traction jump across $\bnd_t$ due to the interior and exterior flows by $\reftracj$. The weak formulation of the direct boundary integral equation becomes
\par
\emph{given $\refconfig{g} \in H^{1/2}(\refconfig{\bnd})$, find $(\reftracj, \zeta) \in H^{-1/2}(\refconfig{\bnd}) \times \IR$ such that}
\begin{equation}
  \label{ch:sto,eq:bieonrefconfig}
    a_{\mathrm{f}t}( \reftracj, \refconfig{\psi} )
  + b_t( \reftracj, \upsilon )
  + b_t( \refconfig{\psi}, \zeta )
  = F_{\mathrm{f}t}( \refconfig{g}, \refconfig{\psi} )
  \qquad \forall (\refconfig{\psi}, \upsilon) \in H^{-1/2}(\bnd_t) \times \IR
  \end{equation}
with bilinear forms
\begin{align*}
  a_{\mathrm{f}t}( \phi, \psi ) &:=
  \langle V_t \phi, \psi \rangle_{H^{1/2}(\refconfig{\bnd})}, \\
  b_t( \psi, \upsilon ) &:=
  \upsilon \langle \psi, J_t n \circ \config_t \rangle_{H^{1/2}(\Gwetref)} +
  \upsilon \langle \psi, J_t n \circ \config_t \rangle_{H^{1/2}(\Gin)}, \\
  F_{\mathrm{f}t}( \phi, \psi ) &:=
  \langle (1+\lambda) J_t \phi/2 + (1-\lambda) K_t \phi, \psi \rangle_{H^{1/2}(\refconfig{\bnd})},
  \end{align*}
and boundary integral operators
\begin{align*}
  \displaystyle
  ( V_t \psi)_i &:=
  J_t \frac{1}{8\pi} \oint\nolimits_\refconfig{\bnd}
  \left( \delta_{ik} \frac{1}{r}  + \frac{(x_i-y_i)(x_k-y_k)}{r^3} \right)
  \circ (\config_t \times \config_t)
  \psi_k(y) J_t \ddmu(y),
  \\
  \displaystyle
  ( K_t \phi)_i &:= 
  J_t \frac{3}{4\pi} \oint\nolimits_\refconfig{\bnd}
  \left( \frac{(x_i-y_i)(x_k-y_k)(x_j-y_j)n_j(y)}{r^5} \right)
  \circ (\config_t \times \config_t)
  \phi_k(y) J_t \ddmu(y),
  \end{align*}
with $i,k\in\{0,1,2\}$, $r=|x-y|$ the Euclidean distance and
\begin{equation}
  \label{ch:sto3,eq:surfjac}
  J_t := (
  (\refconfig{t}^0 \cdot \nabla_{\Gwetref}\config_t) \times
  (\refconfig{t}^1 \cdot \nabla_{\Gwetref}\config_t)
  ) \cdot \refconfig{n}
  \end{equation}
the determinant of the Jacobian of the map $\config_t$. In the above, $\refconfig{n}$ and $\refconfig{t}^\alpha$, $\alpha\in\{0,1\}$ denote the normal and tangent vectors to $\refconfig{\bnd}$. Furthermore, the composition with $(\config_t\times\config_t)$ serves to transport the kernel of the integral operator to the reference manifold.
\par
In~(\ref{ch:sto,eq:bieonrefconfig}), the Lagrange multiplier $\zeta$ constrains the kernel of the boundary-integral operator $V_t$. This kernel is related to the undefined average pressure level $p_0$, similar to the Stokes PDE with pure Dirichlet boundary conditions. Furthermore, $0\leq\lambda\leq\infty$ denotes the viscosity ratio between the fluids interior and exterior to $\bnd$, where the respective limits $\lambda\to0$ and $\lambda\to\infty$ correspond to the separate interior and exterior problems.
\par
One of the advantages that the Stokes model gives is an automatic mechanism for contact prevention, namely, lubrication. It is shown in~\cite[Thms.~33,40]{phdthesis13} that, as two smooth material boundaries advance each other, the fluid opposes this motion with a repulsive traction of singular strength of $O(h^{-3})$, where $h$ is the distance between the boundaries. Thus, the fluid formulation has a built in contact prevention mechanism which avoids the highly nontrivial explicit treatment of self-contact which characterizes folded inflatable structures.
\par
We conclude this section with a number of remarks:
\begin{enumerate}
  \item The boundary integral formulation~(\ref{ch:sto,eq:bieonrefconfig}) gives a direct relation between the Dirichlet and Neumann data at the boundary, mapping velocities imposed by the structure to the tractions imposed on the structure, which is precisely the relation typically required in a fluid-structure interaction problem. No meshing and approximation of the interior domain is thus required, although the solution in the interior can be reconstructed \emph{a posteriori} if desired.
  \item Contrary to a discretization based on a volumetric formulation, it is not more expensive to treat the flow in the exterior domain concurrently, it simply amounts to adjusting the ratio $\lambda$. Moreover, when the conditions in both fluids are the same, $\lambda=1$ and we see that the problem actually becomes cheaper to approximate, as the dual layer operator $K_t$ need not be assembled.
  \end{enumerate}
\par
\subsection{Parametrization-free Kirchhoff-Love structure}\label{sec:structure}
We now turn to the equation governing the configuration $\config$, which is assumed to be a member of a Bochner space. A Bochner space $L^2(0,T;Y)$ contains functions $f:[0,T]\to Y$, such that \smash{$\int_0^T\|f(t,\cdot)\|_Y^2\,dt<\infty$}. First and second time derivatives are denoted $\dot{(\cdot)}$ and $\ddot{(\cdot)}$ respectively. We set $X := H^2(\Gwetref) \cap H_0^1(\Gwetref)$ and write the structure problem as
\par
\emph{given $(\xinit, \vinit) \in H^2(\Gwetref) \times L^2(\Gwetref)$ and $F \in L^2(0,T;L^2(\Gwetref))$, find $\config \in (\xinit + \{\config\in{}L^2(0,T;X):\dot{\config}\in{}L^2(0,T;L^2(\Gwetref)\})$ such that}
\begin{subequations}
  \label{ch:sto3,eq:weakstructure}
  \begin{align} &
  \langle \ddot{\config}_t, \rho \rangle_{X} +
  W'(\config_t; \rho) =
  F(t; \rho)
  \qquad a.e. \; t \in (0,T), \; \forall \rho \in X,
  \\ &
  \config_0 = \xinit,
  \\ &
  \dot{\config}_0 = \vinit.
  \end{align}
  \end{subequations}
In this problem, $F$ is an external load, $\xinit$ and $\vinit$ are the initial position and velocity, and $W$ is the internal energy of a Kirchhoff-Love shell, given in terms of a configuration $\config$ as
\begin{equation}
\label{eq:energy:parfree}
W(\config) = \frac{1}{2} \int_{\Gwetref}
\constit^{ijkl} \left(
\varepsilon_{ij} (\config) \varepsilon_{kl} (\config) +
\upepsilon^2 \kappa_{ij} (\config) \kappa_{kl} (\config)
\right) \ddmu,
\end{equation}
where
\begin{align*}
\varepsilon_{ij} (\config) &:=
\tfrac{1}{2} (
\refconfig{\Pi}_{ij} -
\gammagradi{\config}{m}{i} \gammagradi{\config}{m}{j}), \\
\kappa_{ij} (\config) &:=
\refconfig{\Pi}_{mi} \gammagradi{\refconfig{n}}{m}{j} -
\gammagradi{\config}{m}{i} \gammagradi{n}{m}{j}, \\
\end{align*}
are the components of the membrane and bending strains; and
\begin{equation*}
\constit^{ijkl} := \upsilon
\delta_{ij} \delta_{kl} +
\tfrac{1}{2} (1-\upsilon) (
\delta_{ik} \delta_{jl} +
\delta_{il} \delta_{jk}),
\end{equation*}
is the constitutive tensor of the Saint Venant-Kirchhoff material law, which is especially suited to the anticipated small strains and large rotations~\cite{Bischoff04}. The tensor $\refconfig{\Pi}_{ij} := \delta_{ij} - \refconfig{n}_i \refconfig{n}_j$ denotes the projection into the tangent space of $\Gwetref$ and $\gammagradi{f}{i}{j}$ denotes the $j$\textsuperscript{th} component of the gamma gradient of component $i$ of a vector function $f$.
\par
The structural model introduces two model parameters, $\upepsilon$ and $\upsilon$, which represent the flexural rigidity and Poisson's ratio, respectively.
\par
This form of the Kirchhoff-Love energy has been dubbed \emph{parametrization-free} in~\cite{3d} as it is formulated on the stress-free state $\Gwetref$, without reference to some parametrization of this surface. The main advantage is conceptual: eliminating the superfluous and arbitrary parametrization from the formulation. Implementation aspects for this formulation are treated in~\cite{3d}.
\subsection{Transmission conditions and full problem}\label{sec:coupling}
Standard Dirichlet-Neumann coupling is employed between the fluid and the structure. The kinematic condition imposes a velocity on the fluid at the interface and reads
\begin{subequations} \label{eq:transmission}
  \begin{equation} \label{eq:kinematic}
  \refconfig{g} = \begin{cases}
  \dot{\config}_t & \qquad \text{ on } \Gwetref, \\
  q n & \qquad \text{ on } \Gin,
  \end{cases}
  \end{equation}
with $q:\Gin\to\IR$ the normal inflow velocity. The dynamic condition imposes continuity of tractions at the fluid-structure interface
\begin{equation} \label{eq:kinetic}
  F(t; \cdot) = -\varpi \langle J_t\reftracj, \cdot \rangle.
  \end{equation}
  \end{subequations}
with $J_t$ the determinant of $\config_t$, according to (\ref{ch:sto3,eq:surfjac}). The model parameter $\varpi$ can be interpreted as a \emph{coupling strength} and arises in the nondimensionalization~\cite{3d} by agglomoration of several parameters (such as the Young's modulus) from the dimensional models.
\par
In addition to the coupling conditions, an auxiliary condition is to be satisfied, connected to the incompressibility of the fluid encapsulated by the structure. The configuration provided by the structural equation of motion (\ref{ch:sto3,eq:weakstructure}) has to conserve the volume, i.e., $Q_t(\config_t)=0$ with
\begin{equation}
  \label{eq:compatibility}
  Q_t(\config_t) =
  \frac{1}{3} \left( n, x \right)_{\bnd_t} - \left(
  \frac{1}{3} \left( n, x \right)_\refconfig{\bnd} -
  \int_0^t \int_\Gin q \ddmu \: dt \right),
  \end{equation}
where the first term is the current volume, the second the initial volume and the last the total influx until time $t$. This auxiliary constraint is imposed by the Lagrange multiplier method, and it is proven in~\cite{3d} that this Lagrange multiplier can be seen as the total excess pressure $p_0$ required to uniformly expand the inflatable structure to the correct volume. This constant pressure coincides with the undetermined mode constrained from the fluid subproblem~(\ref{ch:sto,eq:bieonrefconfig}).
\par
The above discussion leads to the aggregated FSI model, composed of the weak form of the interior and exterior fluid boundary integral equations~(\ref{ch:sto,eq:bieonrefconfig}), the Lagrangian equation of motion for the structure~(\ref{ch:sto3,eq:weakstructure}) augmented with the compatibility condition and coupled through the transmission conditions~(\ref{eq:transmission}):
\par
\emph{given $q \in L^2(0,T;L^2(\Gin))$ and $(\xinit, \vinit) \in X \times L^2(\Gwetref)$, find $(\config_t, \reftracj, p_0, \zeta) \in (\xinit + X) \times H^{-1/2}(\refconfig{\bnd}) \times \IR \times \IR$ such that}
\begin{subequations}
\label{ch:sto3,eq:fsiproblem}
  \begin{align} &
    \langle \ddot{\config}_t, \rho \rangle_X
  + W'(\config_t; \rho)
  + p_0 Q_t'(\config_t; \rho)
  + \mu Q_t(\config_t)
  + \varpi \langle J_t\reftracj, \rho \rangle_X
  \nonumber \\ & \qquad
  + a_{\mathrm{f}t}( \reftracj, \refconfig{\psi} )
  + b_t( \reftracj, \upsilon )
  + b_t( \refconfig{\psi}, \zeta )
  - F_{\mathrm{f}t}( \dot{\config}_t, \refconfig{\psi} )
  = F_{\mathrm{f}t}( q(t,\cdot)n, \refconfig{\psi} )
  \nonumber \\ &
  \qquad \qquad
  \qquad a.e. \; t \in (0,T), \; \forall (\rho, \refconfig{\psi}, \mu, \upsilon) \in X \times H^{-1/2}(\refconfig{\bnd}) \times \IR \times \IR,
  \\ &
  \config_0 = \xinit,
  \\ &
  \dot{\config}_0 = \vinit,
  \end{align}
  \end{subequations}
where we extend $q$ and $\dot{\config}_t$ by $0$ on $\Gwetref$ and $\Gin$, respectively.
\section{Approximation and solution}\label{sec:approx}
In treating the discretization and solution of the mathematical model~(\ref{ch:sto3,eq:fsiproblem}), emphasis is put on lubrication and contact mechanics. It was mentioned in~\S\ref{sec:fluid} that the magnitude of the traction forces grow as $O(h^{-3})$ with $h$ the gap size between smooth surfaces in contact. Thus, an isogeometric matching-mesh discretization, which provides the required smoothness and geometrical exactness, is presented in~\S\ref{sec:subdivision}. As the gap size $h$ vanishes, local contributions to the boundary integrals become more singular on account of the $1/r$ and $1/r^2$ dependence of the kernels of $V_t$ and $K_t$. The accurate evaluation of integrals in such circumstances is treated in~\S\ref{sec:bem}. Finally, as the fluid-structure coupling is anticipated to be strong when contact occurs, we outline the partitioned iterative procedure for resolving the coupled system in~\S\ref{sec:partitioned}.
\subsection{Subdivision surface approximation spaces}\label{sec:subdivision}
Two considerations have lead to the choice of approximation space for this problem. Firstly, a $C^1$-continuous basis is desired to have a \emph{rotation free} discretization of the shell structure, and facilitates contact treatment through the lubrication effect. The theoretical results depend on a smooth representation of the fluid boundary. In that case, there is a time at which the gap width is much smaller than the diameter of the contact surface, and the fluid response resembles lubrication. Secondly, an unstructured basis is required, to be able to represent complex geometries, such as folded inflatables; and to perform local refinement where it is required to resolve fine geometrical features and near-singular tractions in near contact.
\par
One of the natural candidates is subdivision surfaces with the Catmull-Clark subdivision scheme~\cite{Catmull78}. Recently, (truncated) hierarchical refinement has been presented for this basis in~\cite{Zore14,Wei15}. One of the freely available CAGD packages in which subdivision surface geometries can be designed is Blender\texttrademark, which has been used in this context.
\par
Catmull-Clark subdivision surfaces can be seen as an extension of cubic b-splines to unstructured grids composed of quadrilaterals. Defining the valence of a vertex as the number of edges terminating there, we observe that interior vertices of a Cartesian mesh are of valency four. An interior vertex of valency other than four is called an \emph{extraordinary vertex} (EV). To extend the definition of the spline basis near an EV, a subdivision mask is introduced, in this case the Catmull-Clark mask, which preserves $C^1$-continuity almost everywhere. As such a mask only provides an implicit definition of the basis functions at specific points, Stam's approach~\cite{Stam98} is used to evaluate the basis functions and their gradients in numerical quadrature. This approach involves $l$ levels of virtual refinement, such that the quadrature point lies in a sub-element away from the EV, and can be evaluated as a regular basis function.
\par
As the domains of the fluid and structure subproblems coincide, it is convenient to use the descretization of the structure for the fluid as well. Thus, this approach can be seen as a \emph{matching mesh} approach, such that any geometrical approximation is avoided. The advantages are twofold. Firstly, errors due to geometrical approximation, which are significant in the BEM~\cite{Scott13}, are now entirely eliminated. One can imagine especially in the pervading near-contact modes, that inaccuracies in the tractions due to geometrical approximation may easily become disastrous, due to the $O(h^{-3})$ relation to the gap size $h$. Secondly, we completely eliminate the volumetric meshing step, which is otherwise a formidable challenge in the light of the complex initial geometries and large deformations characterizing inflatable structures.
\subsection{Boundary-element approximation of the fluid}\label{sec:bem}
To evaluate the singular integrals on product domains in the boundary element problem for the fluid, the regularizing transformations described in~\cite{Sauter11} are used. The convergence rates for standard Gauss quadrature rules are known for this set of transformations and class of integrands, cf.~\cite{Sauter11}. This suggests an adaptive quadrature scheme in which the quadrature order is increased until some heuristic criterion is met. Let $F_{\kappa,q}$ be the approximation of the integral over element $\kappa$ with quadrature order $q$. In this work, the criterion $|F_{\kappa,q-1}-F_{\kappa,q}|<10^{-7}$ is used.
\par
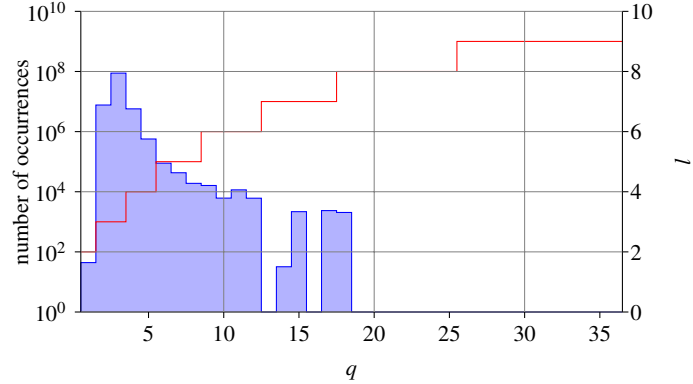
\begin{figure}
  \begin{center}
  \begin{tikzpicture} [scale=0.5, font=\small, x=4mm, y=8mm]
    % plots
    \fill[blue!30, draw=blue] plot file {media/count.plot};
    \draw[red] plot file {media/stam.plot};
    % axis
    \def\xmin{0.5}
    \def\xmax{36.5}
    \def\ymin{0}
    \def\ymax{10}
    \draw (\xmin,\ymin) -- coordinate (x axis mid) (\xmax,\ymin);
    \draw (\xmin,\ymin) -- coordinate (y axis mid) (\xmin,\ymax);
    \draw (\xmax,\ymin) -- coordinate (l axis mid) (\xmax,\ymax);
    \draw[help lines, xstep=10, ystep=2] (\xmin,\ymin) grid (\xmax,\ymax);
    % ticks
    \foreach \x in {5, 10, ..., 35}
      \draw (\x,\ymin) -- (\x,\ymin-0.2)
    node[anchor=north] {\x};
    \foreach \y in {0, 2, ..., 10}
      \draw (\xmin,\y) -- (\xmin-0.4,\y) 
        node[anchor=east] {$10^{\y}$};
    \foreach \y in {0, 2, ..., 10}
      \draw (\xmax+0.4,\y) -- (\xmax,\y) 
        node[anchor=west] {\y};
    % labels
    \node[below=8mm, anchor=center]
      at (x axis mid) {$q$};
    \node[rotate=90, above=8mm, anchor=center]
      at (y axis mid) {number of occurrences};
    \node[rotate=90, below=8mm, anchor=center]
      at (l axis mid) {$l$};
    \end{tikzpicture}
    \end{center}
  \caption{The blue histogram shows the number of occurrences of a quadrature order $q$ from a sample of $10^9$ element evaluations in the simulation of~\S\ref{sec:numex}. The red line shows the required number of levels $l$ in Stam's algorithm to evaluate an order $q$ Gauss rule.}
  \label{fig:stamvsgauss}
  \end{figure}

It is observed that this adaptive procedure converges and that the amount of evaluations requiring an order $q$ decreases rapidly with $q$, cf. fig.~\ref{fig:stamvsgauss}: selection of a scheme of $q>10$ is approximately $10^4$ times less likely than selection of the scheme $q=3$. In fig.~\ref{fig:stamvsgauss} it can also be seen how many levels would be required if an order $q$ scheme would be evaluated near an extraordinary vertex. Such an evaluation would require the $(l-1)$\textsuperscript{th} power of the subdivision matrix, which can be cached once per level. Observe that the number of required levels grows only very slowly and that in this case 8 such levels are required.
\subsection{Partitioned iterative solution}\label{sec:partitioned}
Monolithic solution of a discretization of~(\ref{ch:sto3,eq:fsiproblem}) is untractable, as this requires linearization of the coupling terms. In fact, this leads to hypersingular integrals, i.e., integrands with singularities stronger than those in~(\ref{ch:sto3,eq:fsiproblem}). Despite the strongly coupled nature of this problem due to lubrication, a partitioned scheme is therefore preferred. Thus, the fluid and structure subproblems, the $(\refconfig{\psi},\upsilon)$ and $(\rho,\mu)$ terms of~(\ref{ch:sto3,eq:fsiproblem}) respectively, are solved separately. The required damping for stability of the coupled problem is introduced through the dissipative implicit Euler time-integration scheme applied to the structure.
\par
A subiteration step within our partitioned procedure starts with a linear extrapolation of the initial data (provided by the solution at the previous time interval), which serves as a first approximation of the new coupled solution. Within a fluid--structure subiteration, a structural solve is performed first, to ensure compatibility of the fluid boundary data. The subiteration is considered converged if (1) the norm of the structure residual is below the tolerance before a Newton solve is performed and, in addition, (2) the norm of the fluid update is below that same tolerance.
\section{Deflation of a balloon}\label{sec:numex}
The prime objective of the numerical experiment presented here is to assess the accuracy of the solution method to resolve the complex self-contact that often arises in the simulation of inflatable structures. To this end, as in~\cite{cm12}, an inverted problem is considered, namely the deflation of a balloon.
\par
\begin{figure}
  \centering
  \includegraphics[width=0.8\textwidth]{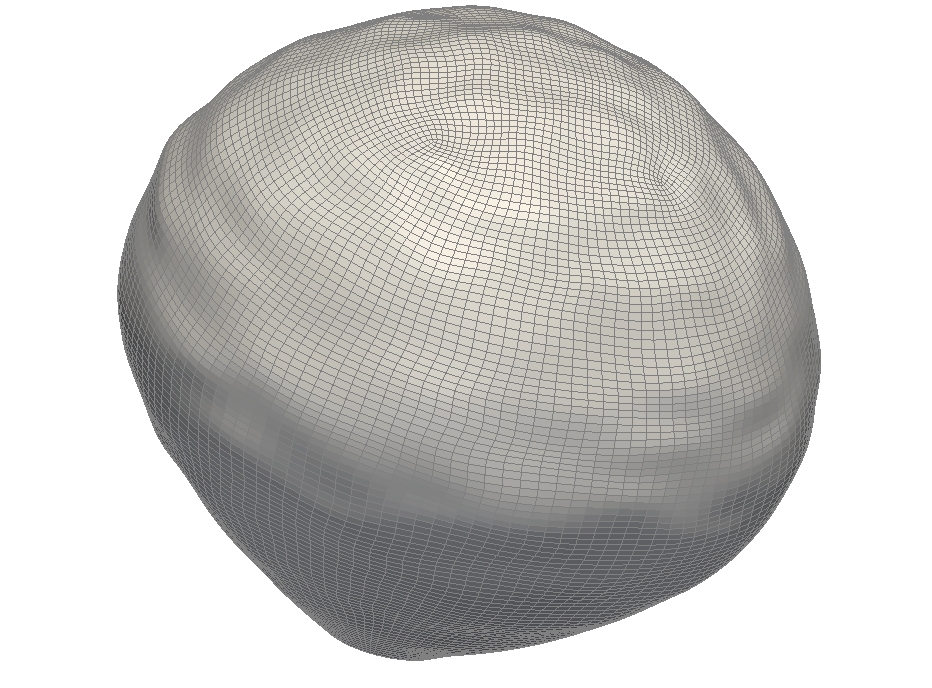}
  \caption{Initial configuration of for the deflation of a balloon.}
  \label{fig:initial}
  \end{figure}

The initial geometry $\refconfig{\bnd}$ contains 320 elements in $\Gin$ and 832 elements in $\Gwetref$. Both these subsets contain 4 extraordinary vertices of valence 3. A small random perturbation is applied to $\Gwetref$ away from $\partial\Gwetref$ so that $\refconfig{\bnd}$ remains connected. The random perturbations causes the structure to wrinkle instead of contract uniformly, rendering the structural response closer to reality. The initial configuration is plotted in fig.~\ref{fig:initial}.
\par
The model parameters are $\upsilon=0$; $\upepsilon=5.77\;10^{-4}$ corresponding to a flexible shell; $\lambda=1$ corresponding to identical fluids in the interior and exterior; and $\varpi=10^{-5}$. The outflow $q$ is constant in time and has a sine-shaped profile with a magnitude such that the volume would vanish at $T=2^{12}$. The time step size is $\tau=4$ and the tolerance in the partitioned solver \texttt{TOL} is set to $10^{-6}$, settings for which only 1-2 subiterations are required at each time level.
\begin{figure}
  \begin{center}
  \subfloat[$t=4.0\;10^2$, $m(\dom_t)/m(\dom_0)=0.80$]{ \label{fig:first}
    \snapshot{$\reftracj_2$}{$-3$}{$3$}{ 
      \includegraphics[width=0.45\textwidth]{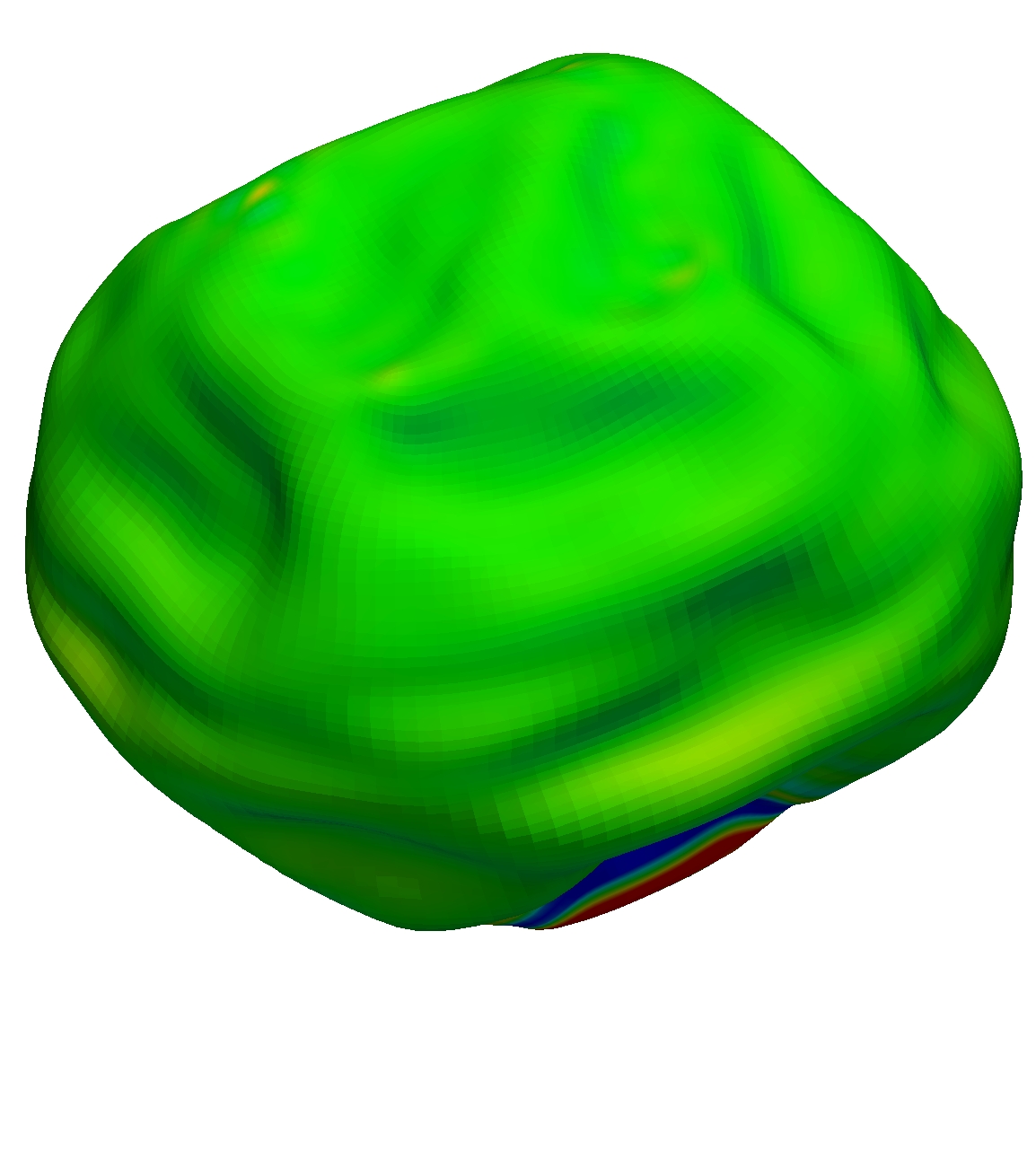} }}
  \subfloat[$t=8.0\;10^2$, $m(\dom_t)/m(\dom_0)=0.61$]{
    \snapshot{$\reftracj_2$}{$-3$}{$3$}{
      \includegraphics[width=0.45\textwidth]{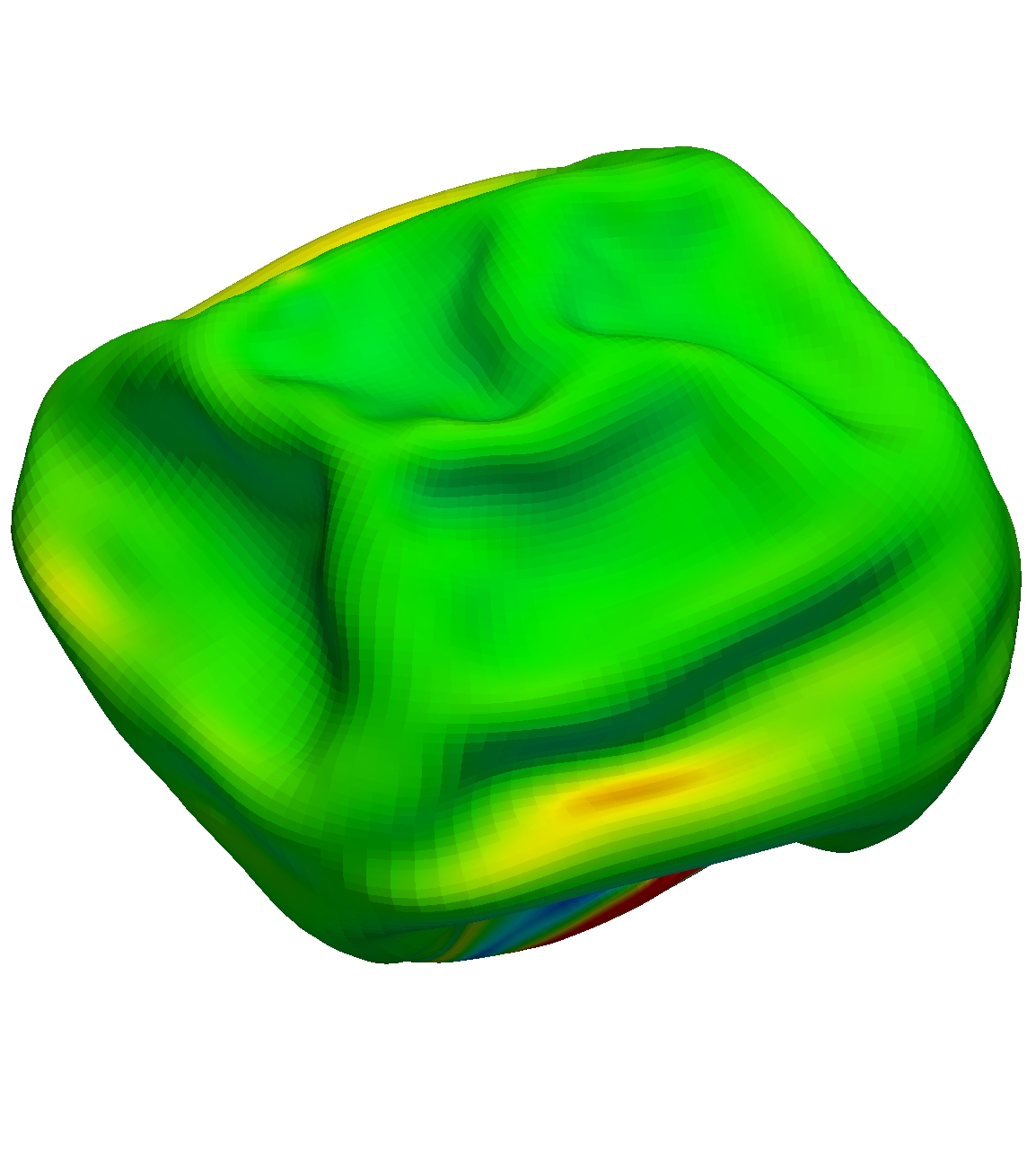} }} \\
  \subfloat[$t=1.20\;10^3$, $m(\dom_t)/m(\dom_0)=0.41$]{ \label{fig:penult}
    \snapshot{$\reftracj_2$}{$-3$}{$3$}{
      \includegraphics[width=0.45\textwidth]{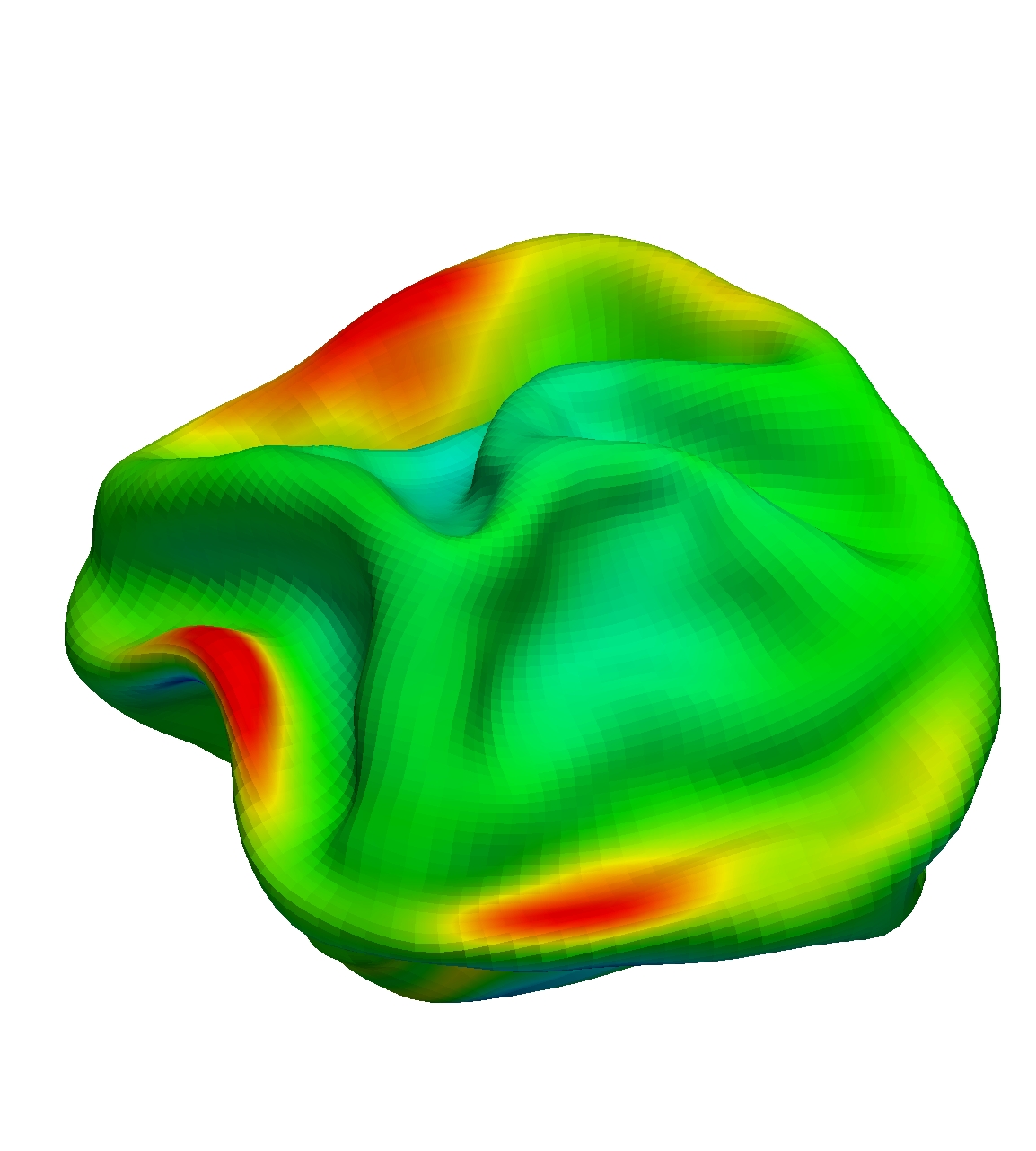} }}
  \subfloat[$t=1.31\;10^3$, $m(\dom_t)/m(\dom_0)=0.26$]{ \label{fig:last}
    \snapshot{$\reftracj_2$}{$-3$}{$3$}{ 
      \includegraphics[width=0.45\textwidth]{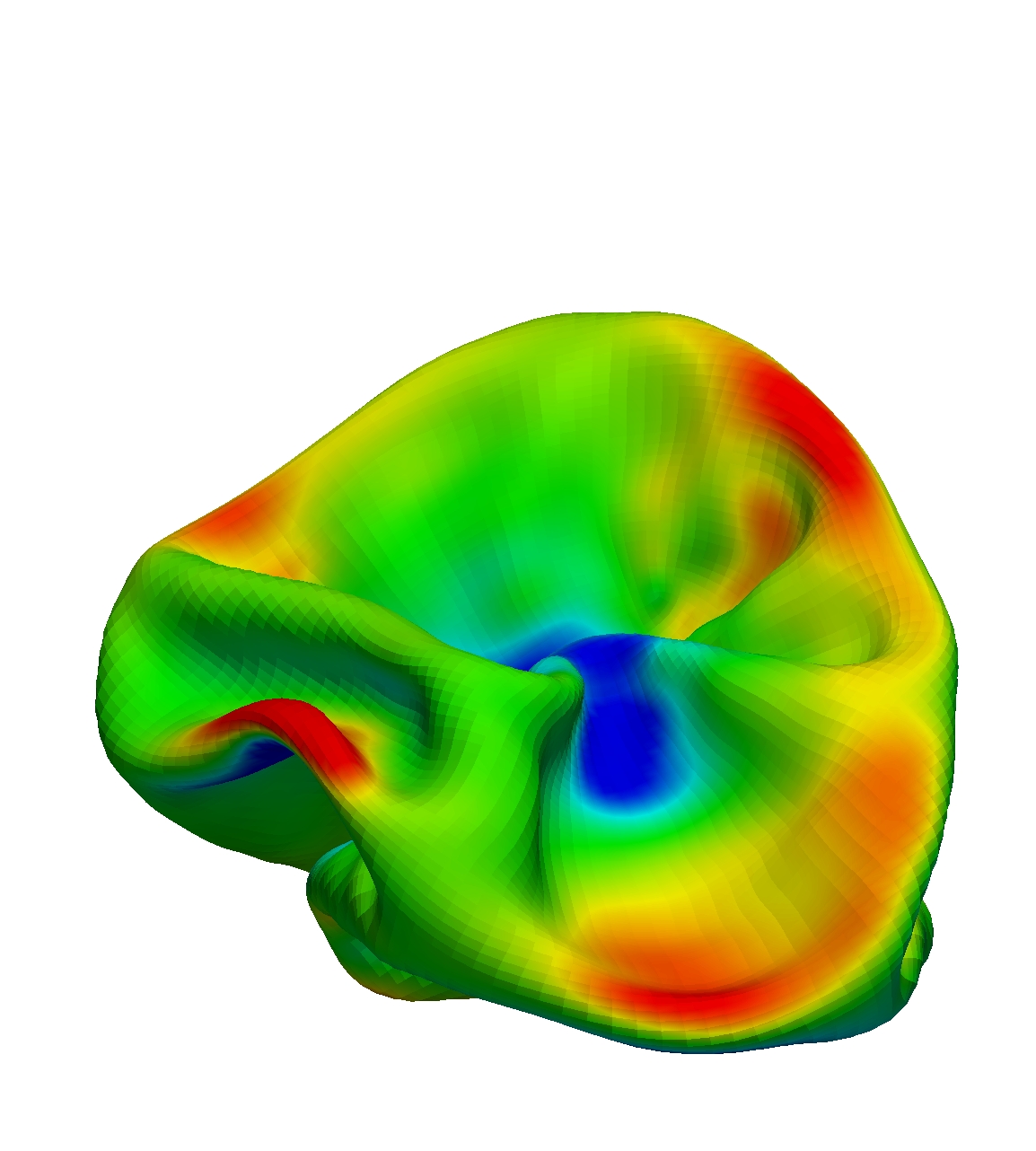} }} \\
  \end{center}
  \caption{Snapshots of the configuration and fluid traction at different time levels during the deflation process. The color coding corresponds to the vertical component of the traction, $\reftracj_2$ and $m(A)$ denotes the volume of set $A$.}
  \label{fig:snapshots}
  \end{figure}
\par
In fig.~\ref{fig:snapshots}, snapshots of the deflation process are shown. Recall that the plotted traction component $\reftracj_2$ does not include the contribution of the total excess pressure $p_0n$. The action of the lubrication effect can clearly be observed along horizontal folds, where the vertical component is positive (red) on upper surfaces and negative (blue) on lower surfaces and close to zero elsewhere in figs.~\ref{fig:first}--\ref{fig:penult}. Thus, the collapse of the structure under the action of the total excess pressure $p_0n$, which is uniform over $\Gwet$, is observed to be counteracted locally by $\reftracj$ to prevent contact.
\par
As the limits of the color bars are identical, it is clear that the magnitude of these lubrication forces increases as the deflation proceeds. At the same time, the folds in the fabric become gradually sharper. The deflation process is terminated at $t=1.31\;10^3$, cf. fig.~\ref{fig:last}. As can be seen from this snapshot, a large downward traction is exerted on the part of the structure near $\Gin$, where the outflow condition~(\ref{eq:kinematic}) is prescribed. The outflux $q$ does not correspond to the velocity of $\Gin$. The segments therefore seem to recede and no lubrication force is instigated. Ultimately, this self-intersection is therefore not due to insufficient resolution of the fluid response, but rather to the artificial problem setup where a structure is deflated instead of inflated.
\section{Conclusions}\label{sec:concl}
In this chapter we have presented a FEBE method for the simulation of inflatable structures. A Kirchhoff-Love shell with low flexural rigidity, discretized with the finite element method (using a Catmull-Clark subdivision surface basis) is coupled with a boundary element method discretization of both the interior and exterior fluid. It is advantageous to let the fluid inherit the structural mesh, as this leads to a geometrically exact matched discretization at the fluid-structure interface, which in turn enables accurate approximation of contact forces. The coupled system may be advanced in a partitioned iterative fashion.
\par
We have furthermore assessed the capabilities of this method in the presence of multiple modes of self-contact, by considering a deflation problem, in which the fluid response is anticipated to be dominated by lubrication forces. It was observed that the computed traction forces effectively prevent self-intersection of the structure, although very large deformations are sustained.
\bibliographystyle{plain} % if abbrv etc used, Dutch surnames are incorrect
{\renewcommand\newpage{}\bibliography{o}}

\begin{thebibliography}{10}

\bibitem{Bischoff04}
M.~Bischoff, W.~A. Wall, K.-U. Bletzinger, and E.~Ramm.
\newblock {\em Encyclopedia of Computational Mechanics}, volume 2 Structures,
  chapter 3 Models and Finite Elements for Thin-walled Structures, pages
  59--137.
\newblock Wiley, 2004.

\bibitem{Catmull78}
E.~Catmull and J.~Clark.
\newblock Recursively generated b-spline surfaces on arbitrary topological
  meshes.
\newblock {\em Computer-Aided Design}, 10(6):350--355, 1978.

\bibitem{Hsu14a}
M.-C. Hsu, D.~Kamensky, Y.~Bazilevs, M.S. Sacks, and T.J.R. Hughes.
\newblock Fluid--structure interaction analysis of bioprosthetic heart valves:
  significance of arterial wall deformation.
\newblock {\em Computational Mechanics}, 54:1055--1071, 2014.

\bibitem{Hughes81}
T.J.R. Hughes, W.K. Liu, and T.K. Zimmermann.
\newblock {L}agrangian-{E}ulerian finite element formulation for incompressible
  viscous flows.
\newblock {\em Computer Methods in Applied Mechanics and Engineering},
  29:329--349, 1981.

\bibitem{Johnson94}
A.A. Johnson and T.E. Tezduyar.
\newblock Mesh update strategies in parallel finite element computations of
  flow problems with moving boundaries and interfaces.
\newblock {\em Computer Methods in Applied Mechanics and Engineering},
  119:73--94, 1994.

\bibitem{Johnson99}
A.A. Johnson and T.E. Tezduyar.
\newblock Advanced mesh generation and update methods for 3d flow simulations.
\newblock {\em Comput. Mech.}, 23:130--143, 1999.

\bibitem{Kamensky15}
D.~Kamensky, M.-C. Hsu, D.~Schillinger, J.A. Evans, A.~Aggarwal, Y.~Bazilevs,
  M.S. Sacks, and T.J.R. Hughes.
\newblock An immersogeometric variational framework for fluid--structure
  interaction: Application to bioprosthetic heart valves.
\newblock {\em Computer Methods in Applied Mechanics and Engineering},
  284:1005--1053, 2015.

\bibitem{phdthesis13}
{T.M. van} Opstal.
\newblock {\em Numerical methods for inflatables with multiscale geometries}.
\newblock PhD thesis, Eindhoven University of Technology, 2013.

\bibitem{cmame12}
{T.M. van} Opstal and {E.H. van} Brummelen.
\newblock A finite-element/boundary-element method for large-displacement
  fluid-structure interaction with potential flow.
\newblock {\em Computer Methods in Applied Mechanics and Engineering},
  266:57--69, http://dx.doi.org/10.1016/j.cma.2013.07.009 2013.

\bibitem{cm12}
{T.M. van} Opstal, {E.H. van} Brummelen, {R. de} Borst, and M.R. Lewis.
\newblock A finite-element/boundary-element method for large-displacement
  fluid-structure interaction.
\newblock {\em Computational Mechanics}, 50(6):779--788, 2012.

\bibitem{3d}
{T.M. van} Opstal, {E.H. van} Brummelen, and {G.J. van} Zwieten.
\newblock A finite-element/boundary-element method for three-dimensional,
  large-displacement fluid-structure interaction.
\newblock {\em Computer Methods in Applied Mechanics and Engineering},
  284:637--663, 2015.

\bibitem{Sauter11}
S.A. Sauter and C.~Schwab.
\newblock {\em Boundary Element Methods}.
\newblock Springer Series in Computational Mechanics. Springer-Verlag, 2011.

\bibitem{Scott13}
M.A. Scott, R.N. Simpson, J.A. Evans, S.~Lipton, S.P.A. Bordas, T.J.R. Hughes,
  and T.W. Sederberg.
\newblock Isogeometric boundary element analysis using unstructured
  {T}-splines.
\newblock {\em Computer Methods in Applied Mechanics and Engineering},
  254:197--221, 2013.

\bibitem{Stam98}
J.~Stam.
\newblock Exact evaluation of {C}atmull-{C}lark subdivision surfaces at
  arbitrary parameter values.
\newblock In {\em proceedings of the 25\textsuperscript{th} annual conference
  on computer graphics and interactive techniques}, pages 295--404. ACM, 1998.

\bibitem{Takizawa12b}
K.~Takizawa, M.~Fritze, D.~Montes, T.~Spielman, and T.E. Tezduyar.
\newblock Fluid-structure interaction modeling of ringsail parachutes with
  disreefing and modified geometric porosity.
\newblock {\em Computational Mechanics}, 50:835--854, 2012.

\bibitem{Takizawa11}
K.~Takizawa, T.~Spielman, and T.E. Tezduyar.
\newblock Space--time {FSI} modeling and dynamical analysis of spacecraft
  parachutes and parachute clusters.
\newblock {\em Computational Mechanics}, 48:345--364, 2011.

\bibitem{Takizawa15}
K.~Takizawa, T.E. Tezduyar, C.~Boswell, R.~Kolesar, and K.~Montel.
\newblock {FSI} modeling of the reefed stages and disreefing of the {O}rion
  spacecraft parachutes.
\newblock {\em Computational Mechanics}, 54:1203--1220, 2014.

\bibitem{Tezduyar92}
T.E. Tezduyar.
\newblock Stabilized finite element formulations for incompressible flow
  computations.
\newblock {\em Advances in Applied Mechanics}, 38:1--44, 1992.

\bibitem{Tezduyar10}
T.E. Tezduyar, K.~Takizawa, C.~Moorman, S.~Wright, and J.~Christopher.
\newblock Space-time finite element computation of complex fluid-structure
  interactions.
\newblock {\em International Journal for Numerical Methods in Fluids},
  64:1201--1218, 2010.

\bibitem{Wei15}
X.~Wei, Y.~Zhang, T.J.R. Hughes, and M.A. Scott.
\newblock Truncated hierarchical catmull--clark subdivision with local
  refinement.
\newblock {\em Computer Methods in Applied Mechanics and Engineering},
  291:1--20, 2015.

\bibitem{Wick13b}
T.~Wick.
\newblock Flapping and contact {FSI} computations with the fluid--solid
  interface-tracking/interface-capturing technique and mesh adaptivity.
\newblock {\em Computational Mechanics}, page published online, 2013.

\bibitem{Zore14}
U.~Zore, B.~J{\"u}ttler, and J.~Kosinka.
\newblock On the linear independence of (truncated) hierarchical subdivision
  splines.
\newblock {\em Geometry+ Simulation Report}, 17, 2014.

\end{thebibliography}
\end{document}